\documentclass[12pt]{amsart}

\usepackage{environ, xcolor}
\NewEnviron{commentA}{}

\newcommand\Aoff{\RenewEnviron{commentA}{}}
\Aoff{} 

\usepackage{amsmath, amssymb}
\usepackage{array}
\usepackage[frame,cmtip,arrow,matrix,line,graph,curve]{xy}
\usepackage{graphpap, color, paralist, pstricks}
\usepackage[mathscr]{eucal}
\usepackage[pdftex]{graphicx}
\usepackage[pdftex,colorlinks,backref=page,citecolor=blue]{hyperref}
\usepackage{cleveref}
\usepackage{verbatim}
\usepackage{makecell}

\newtheorem{theorem}{Theorem}[section]
\newtheorem{proposition}[theorem]{Proposition}

\newtheorem{lemma}[theorem]{Lemma}

\theoremstyle{definition}

\newtheorem{example}[theorem]{Example}

\newtheorem{question}[theorem]{Question}
\newtheorem{remark}[theorem]{Remark}

\usepackage[margin=1in]{geometry} 
\usepackage{amsmath,amssymb,mathtools, mathrsfs,esint,verbatim}
\newcommand{\R}{{\mathbb R}}

\newcommand{\Z}{{\mathbb Z}}
\newcommand{\Q}{{\mathbb Q}}

\newcommand{\C}{{\mathbb C}}

\newcommand{\on}[1]{\operatorname{#1}}

\newcommand{\Spec}{{\on{Spec}}}

\DeclareFontFamily{U}{mathb}{\hyphenchar\font45}
\DeclareFontShape{U}{mathb}{m}{n}{
      <5> <6> <7> <8> <9> <10> gen * mathb
      <10.95> mathb10 <12> <14.4> <17.28> <20.74> <24.88> mathb12
      }{}
\DeclareSymbolFont{mathb}{U}{mathb}{m}{n}
\DeclareMathSymbol{\righttoleftarrow}{3}{mathb}{"FD}


\title{Rational Weighted Projective Hypersurfaces}
\author{Louis Esser}

\address{Department of Mathematics, Princeton University, Fine Hall, Washington Road, Princeton, NJ 08544-1000, USA}
\email{esserl@math.princeton.edu}

\begin{document}

\begin{abstract}
A very general hypersurface of dimension $n$ and degree $d$ in 
complex projective space is rational if $d \leq 2$, but is expected to be 
irrational for all $n, d \geq 3$.  Hypersurfaces in 
weighted projective space with degree small relative to the
weights are likewise rational.
In this paper, we introduce rationality constructions for 
weighted hypersurfaces of higher degree that provide many new rational examples
over any field.  We answer in the affirmative a question of T. Okada
about the existence
of very general terminal Fano rational weighted
hypersurfaces in all dimensions $n \geq 6$.
\end{abstract}

\maketitle

\section{Introduction}

An irreducible algebraic variety $X$ over the field $k$ is 
($k$-)\textit{rational} if there
is a $k$-birational map $X \dashrightarrow \mathbb{P}^{\dim(X)}_{k}$.  
Determining 
whether a given variety is rational is often extremely difficult.  
One classical problem is to determine which hypersurfaces 
in $\mathbb{P}^{n+1}_k$ are rational.

If $X \subset \mathbb{P}^{n+1}_k$ is a hypersurface of degree $1$, it is 
a hyperplane, so it is isomorphic to $\mathbb{P}^n_k$, hence rational.  When $d = 2$,
$X$ is a quadric, which is
rational if and only if it contains a smooth $k$-point.  For $d \geq 3$,
the problem is much subtler. For instance, when $k = \C$, 
it's expected that the very general
hypersurface of degree $d \geq 3$ in $\mathbb{P}^{n+1}_{\C}$ is not rational
for any $n$, with the
exception of cubic surfaces.  However, this is unproven already for cubic 
hypersurfaces in $\mathbb{P}_{\C}^5$.

Recently, the more general question of rationality for hypersurfaces in 
weighted projective space $\mathbb{P}(a_0,\ldots,a_{n+1})$ has also been 
considered \cite{CP,CTP,Okada2,Okada3,Okada}.
Weighted projective hypersurfaces are a diverse class of 
varieties that include hypersurfaces in $\mathbb{P}^{n+1}$ 
as well as natural 
geometric constructions such as cyclic covers of projective space.
In this paper, we find many new examples of rational weighted projective 
hypersurfaces
over any field $k$ by introducing two main rationality
constructions.

The first construction, \Cref{thm:Delsarte},
generalizes a result of J. Koll\'{a}r showing
that certain hypersurfaces in weighted projective space with ``loop" 
equations are rational, among other remarkable properties 
\cite[Section 5]{Kollar}.  The generalization pertains to any hypersurface
defined by a \textit{Delsarte equation}, i.e., an equation with the
same number of monomials as variables.  The idea of the proof is to
construct a birational map to $\mathbb{P}^n$ directly, using the linear
system generated by the monomials.  This works under a certain
gcd condition on the exponents in the equation.

The second construction, \Cref{thm:quadric_bundle}, shows that certain
weighted projective hypersurfaces admit a birational quadric bundle 
structure with a section over a rational variety, and hence are rational.

Just like in ordinary projective space, 
weighted hypersurfaces of
``low degree" are easily shown to be rational. More precisely, a 
quasismooth weighted projective hypersurface 
$X_d \subset \mathbb{P}_{\C}(a_0,\ldots,a_{n+1})$ (with $a_0 \geq a_1 \geq
\cdots \geq a_{n+1}$) is always rational whenever the following criterion
holds (see \Cref{prop:lowdeg} for a more general statement):
\begin{equation}
\label{eq:degree_criterion} \tag{I}
d < 2a_0 \text{ or } d = 2a_0 = 2a_1.
\end{equation}
We use the constructions above to show that, contrary to the expectation 
in ordinary
projective space, there \textit{are} many examples with 
$d > 2\max\{a_0,\ldots,a_{n+1}\}$ where $X$ is very general, 
quasismooth, and rational.
In particular, we answer a question of
T. Okada that arose from his study of rationality for terminal Fano
threefold hypersurfaces \cite{Okada}:

\begin{question}{{{{\cite[Question 1.3]{Okada}}}}}
\label{question:main}
Let $X_d \subset \mathbb{P}_{\C}(a_0,\ldots,a_{n+1})$ be a very general 
quasismooth weighted projective hypersurface of 
dimension 
$n \geq 3$ with terminal singularities.
Can $X$ be rational
without satisfying the degree criterion 
\eqref{eq:degree_criterion}?
\end{question}

We note that \Cref{question:main} is formulated for $n \geq 3$ because
of the example of cubic surfaces in $\mathbb{P}_{\C}^3$,
which are rational but fail
the criterion \eqref{eq:degree_criterion}.  In addition, there are 
several other families of singular quasismooth
weighted projective surfaces which are 
known to be rational, some of which even have ample canonical divisor
(see \cite{Chitayat}).

The answer to \Cref{question:main} is negative when 
$n = 3$, meaning that the degree criterion is necessary and 
sufficient for rationality in that case.
Indeed, there are a total of 
$130$ families of quasismooth terminal weighted
projective Fano threefold hypersurfaces, among which $20$ satisfy
\eqref{eq:degree_criterion}, and hence are 
rational.  Of the remaining $110$ families, $95$ have Fano index $1$. 
I. Cheltsov and J. Park proved that every quasismooth member of each of these
$95$ families is birationally rigid, and in particular not rational \cite{CP}.
Okada showed that a very general member of the other $15$ families is not 
rational \cite[Theorem 1.1]{Okada}. In fact, with the lone exception of 
cubic hypersurfaces in $\mathbb{P}_{\mathbb{C}}^4$, he proved that
a very general hypersurface from any of the $110$ families failing
\eqref{eq:degree_criterion} is not stably rational.

We show that the answer to 
\Cref{question:main} is actually affirmative in all dimensions $n \geq 6$
(and for all $n \geq 3$ if we weaken ``terminal" to ``klt").

\begin{theorem}
\label{thm:main}
For every integer
$n \geq 3$, there exist positive integers 
$d,a_0,\ldots,a_{n+1}$ such that $d > 2 \max\{a_0,\ldots,a_{n+1}\}$
and every quasismooth
$X_d \subset \mathbb{P}_{\C}(a_0,\ldots,a_{n+1})$ is 
a rational klt Fano variety
(and there exist such quasismooth $X$).
For $n \geq 6$, we can also
make $X$ terminal.
\end{theorem}

The terminal cases for $n = 4,5$ remain open.
Though the question was originally formulated over the complex numbers,
our methods yield quasismooth rational examples 
with $d > 2 \max\{a_0,\ldots,a_{n+1}\}$
in every dimension over any field $k$.
We also find examples with non-trivial moduli.
We note that there is another approach to answering
\Cref{question:main}, at least for $n \geq 7$,
using a different rationality construction
due to M. Artebani and M. Chitayat (see \Cref{rem:alternative}).
It's interesting to note that the various constructions
all seem to produce terminal examples beginning only in
dimension $6$ or $7$.

\Cref{sect:prelim} defines the necessary terminology related to 
weighted projective varieties.  \Cref{sect:ratl} establishes
rationality criteria for weighted projective hypersurfaces.
In \Cref{sect:ex}, we find some special choices for
$d,a_0, \ldots,a_{n+1}$ that
produce examples proving \Cref{thm:main}.  We also give
a method of checking whether loop hypersurfaces are canonical
or terminal using only the exponents.
Finally, we show that there are families with non-trivial moduli
satisfying the conditions of \Cref{thm:main}.

\noindent \textit{Acknowledgements.} We thank J\'{a}nos Koll\'{a}r and Burt Totaro
for many helpful suggestions and comments.

\section{Preliminaries}
\label{sect:prelim}

Over a base field $k$, we define the \textit{weighted projective space} with 
weights $a_0,\ldots,a_{n+1}$ to be the projective variety 
$\mathbb{P} \coloneqq \mathbb{P}^{n+1}_k(a_0,\ldots,a_{n+1})
\coloneqq \mathrm{Proj} (k[x_0,\ldots,x_{n+1}])$.  Here the $a_i$
are positive integers and the variable $x_i$ has weight $a_i$.
We sometimes use the abbreviation $a^{(r)}$ to indicate 
that the weight $a$ appears $r$ times.
The grading of the polynomial ring above
corresponds to a $\mathbb{G}_{\on{m},k}$-action on affine space
$\mathbb{A}_k^{n+2} = \Spec(k[x_0,\ldots,x_{n+1}])$ given by the ring 
homomorphism
\begin{align*}
    k[x_0,\ldots,x_{n+1}] & \longrightarrow  k[x_0,\ldots,x_{n+1}] \otimes_k 
    k[t,t^{-1}], \\
    x_i & \longmapsto x_i \otimes t^{a_i}.
\end{align*}
Then $\mathbb{P}_k(a_0,\ldots,a_{n+1})$ is also identified with the 
universal quotient $(\mathbb{A}_k^{n+2} \setminus \{0\}) / \mathbb{G}_{\on{m},k}$.
There is a third description of $\mathbb{P}$ as a toric variety: it is the toric
variety over $k$
associated to the fan consisting of cones generated by proper subsets
of the vectors $e_0, \ldots, e_n, v_{n+1} \in \mathbb{R}^{n+1}$ 
in the lattice $N$ generated by these vectors, where
$$v_{n+1} \coloneqq -\frac{a_0}{a_{n+1}}e_0 - \cdots - \frac{a_n}{a_{n+1}}e_n.$$ 
We call a set $\{c_1,\ldots,c_r\}$ of integers \textit{well-formed} if 
$\gcd(c_1,\ldots,\widehat{c_i},\ldots,c_r) = 1$ for each $i = 1,\ldots,r$.  A
weighted projective space is \textit{well-formed}
if its set of weights is so.  Any
weighted projective space is isomorphic, as a variety, to one which is 
well-formed.

Via the Proj construction, $\mathbb{P}$ is equipped with a reflexive sheaf 
$\mathcal{O}(1)$, which is associated to a Weil divisor on $X$. 
In the well-formed case, the canonical class of $\mathbb{P}$ is
$K_{\mathbb{P}} = \mathcal{O}(-\sum_i a_i)$.

A subvariety $X \subset \mathbb{P}$ is \textit{quasismooth} if its preimage
$C_X^* \subset \mathbb{A}_k^{n+2} \setminus \{0\}$ under the quotient 
morphism is 
smooth over $k$.  The affine variety $C_X^*$ is called the 
\textit{(punctured) affine cone} over $X$.
The subvariety $X$ is \textit{well-formed} if $\mathbb{P}$ is well-formed 
and the intersection of $X$ with the 
singular locus of $\mathbb{P}$ has codimension at least $2$ in $X$. 
A \textit{weighted projective hypersurface} is a subvariety of $\mathbb{P}$
defined by a single polynomial $f(x_0,\ldots,x_{n+1})$
with $k$-coefficients that is 
weighted homogeneous of degree $d$.
A well-formed quasismooth hypersurface of degree $d$ satisfies the 
adjunction formula
$K_X = \mathcal{O}_X(d - \sum_i a_i)$.

A few special types of polynomials $f$ appear in the paper.  A polynomial $f$
is \textit{Delsarte} if it has the same number of monomials as variables.
One example is a \textit{loop polynomial}, which has the form:
$$x_0^{b_0} x_1 + x_1^{b_1} x_2 + \cdots + x_{n+1}^{b_{n+1}}x_0.$$
We say that this is a loop polynomial of \textit{type} 
$[b_0,\ldots,b_{n+1}]$.

Unless all the weights equal $1$, a well-formed weighted projective
space is always singular, but its singularities are of a special class called
cyclic quotient singularities. A \textit{cyclic quotient singularity
of type} $\frac{1}{r}(c_1,\ldots,c_s)$ is a singular point 
\'{e}tale-locally isomorphic to the point $0 \in \mathbb{A}^s/\mu_r$, where the group $\mu_r$ acts by 
$\zeta \cdot (t_1,\ldots,t_s) = (\zeta^{c_1}t_1,\ldots,\zeta^{c_s}t_s)$
for any $\zeta \in \mu_r$.  We say this singularity is \textit{well-formed}
if 
$\gcd(r,c_1,\ldots,\widehat{c_j},\ldots,c_s) = 1$ for all $j = 1,\ldots,s$.

For the rest of this section take $k = \C$. Then a quasismooth weighted
projective hypersurface also has only cyclic quotient singularities, whose
types are determined by the weights and degree of the hypersurface (see 
\cite[Lemma 2.5, Proposition 2.6]{ETW}).  There is a combinatorial 
condition called the \textit{Reid-Tai criterion}
that can be used to determine whether a cyclic quotient singularity of
a particular type belongs to certain classes important
to the Minimal Model Program.  In the theorem below, we use $[x]$ to 
denote $x - \lfloor x \rfloor$, the fractional part of $x$.

\begin{theorem}{{{\cite[Theorem 4.11]{Reidyoung}}}}
\label{thm:RT}
Let $\frac{1}{r}(c_1,\ldots,c_s)$ be a well-formed cyclic quotient singularity. 
This singularity is canonical (resp.\ terminal) if and only if
$$\sum_{j = 1}^s \left[\frac{i c_j}{r} \right] \geq 1 $$
(resp.\ $> 1$) for all $i = 1,\ldots,r-1$.
\end{theorem}

We also note that all quotient singularities (over $\C$) are klt.

\section{Rationality criteria for weighted projective hypersurfaces}
\label{sect:ratl}

Since every weighted projective space is a toric variety, it is rational 
(over any field). 
Hypersurfaces for which the degree is small compared to the weights
are also rational.  The proposition below is a generalization of 
the degree criterion \eqref{eq:degree_criterion} of T. Okada
to an arbitrary field $k$.

\begin{proposition}
\label{prop:lowdeg}
Suppose that $X_d \subset \mathbb{P}_k(a_0,\ldots,a_{n+1})$ is a well-formed 
quasismooth hypersurface over any field $k$.  Assume one of
the following two conditions holds:
\begin{enumerate}
    \item \label{prop:lowdeg-1} $d < 2a_0$;
    \item \label{prop:lowdeg-2} $d = 2a_0 = \cdots = 2a_r$,
    for some $r \geq 1$, and $X$ contains a point in 
    $\{x_{r+1} = \cdots = x_{n+1} = 0\}$ defined over $k$.
\end{enumerate}
Then $X$ is rational over $k$.
\end{proposition}

In particular, a quasismooth 
$X_d \subset \mathbb{P}_k(a_0,\ldots,a_{n+1})$ is always
rational over an algebraically closed field $k$ when $d \leq 2
\max\{a_0,\ldots,a_{n+1}\}$, unless $d = 2a_0$ and only one weight equals 
half the degree. There are many non-rational examples in the latter case, 
such as double covers $X_{2a} \subset \mathbb{P}(a,1^{(n+1)})$ of 
$\mathbb{P}^n$ branched in a divisor of degree $2a$ with $a \geq n+1$.
(Indeed, $X$ is smooth and $K_X = \mathcal{O}_X(c)$ has global sections 
in this
case since $c \coloneqq 2a - a - n - 1 \geq 0$.)

\begin{proof}
Let $f$ be the weighted homogeneous polynomial defining the hypersurface $X$.
First suppose that \eqref{prop:lowdeg-1} holds.  If there is any monomial of 
the form $x_i$ in $f$, then we may write the equation $f = 0$ as $x_i = 
g(x_0,\ldots,\widehat{x_i},\ldots,x_{n+1})$.  In this case, $X$ is called a 
\textit{linear cone}, and the change of coordinates $x_i \mapsto x_i - g$, 
defined over $k$, exhibits an isomorphism from $X$ to $\{x_i = 0\} \cong 
\mathbb{P}(a_0,\ldots,\widehat{a_i},\ldots,a_{n+1})$, which is rational.

If there is no monomial of the form $x_i$ in $f$ and $d < 2a_0$ holds,
then there must 
be a monomial of the form $x_0 x_i$ in $f$, or else the affine cone over $X$ 
would contain $(1,0,\ldots,0) \in \mathbb{A}^{n+2}_k$ and be singular at this 
point, contradicting quasismoothness.  After an appropriate coordinate change 
in $x_i$, defined over $k$, we can assume that $x_0 x_i$ is the only
monomial involving $x_0$, so that $f = 0$ can be written
$$x_0 x_i = g(x_1,\ldots,x_{n+1}).$$
\begin{commentA}
In detail, we begin with an equation of the form $x_0 x_i = x_0 g_1 + g_2$,
where $g_1$ doesn't involve $x_0$ or $x_i$, and $g_2$ does not involve $x_0$.
We set $x_i \mapsto
x_i - g_1$ to eliminate 
the $g_1$ term.  This puts the equation in the required form.  
Since $g_1$ has
$k$-coefficients, this coordinate changes is defined over $k$.

\end{commentA}
There is then a $k$-birational map 
$X \dashrightarrow \mathbb{P}(a_1,\ldots,a_{n+1})$ defined
by forgetting
the coordinate $x_0$.  Hence $X$ is rational.

Now assume \eqref{prop:lowdeg-2} holds.  Then $f$ is a sum of 
a quadratic form in $x_0, \ldots, x_r$ and terms involving other variables.
The closed stratum where $x_{r+1},\ldots,x_{n+1}$ vanish is isomorphic to
$\mathbb{P}^r$, and by assumption its intersection with $X$ contains a 
$k$-point.  After a change of variables in $x_0,\ldots,x_r$,
we may assume this $k$-point has $x_0 \neq 0$ but $x_1 = \cdots = x_r = 0$.

Since $X$ contains this point, there is no monomial of the form $x_0^2$ in $f$.
The quasismoothness condition at this point therefore requires that there is a monomial
of the form $x_0x_i$ with nonzero coefficient in $f$.  From here, we can
use the same argument from the proof of \eqref{prop:lowdeg-1} to show that
$X$ is rational.
\end{proof}

The remainder of this section presents two new rationality criteria
for weighted projective hypersurfaces.  These in particular produce
many quasismooth rational examples with $d > 2 \max\{a_0,\ldots,a_{n+1}\}$.

\subsection{Delsarte hypersurfaces}
\label{subsect:Delsarte}

A \textit{Delsarte} polynomial is a weighted homogeneous polynomial
with the same number of monomials as variables.  Given a Delsarte 
polynomial $f$ defining a weighted projective hypersurface in 
$\mathbb{P}(a_0,\ldots,a_{n+1})$, we can associate to it an 
$(n+2) \times (n+2)$ matrix 
$B = (b_{ij})$, where the entries $b_{ij}$ are determined 
from the equation $f$ as follows:
$$f = \sum_{i = 0}^{n+1} K_i \prod_{j = 0}^{n+1} x_j^{b_{ij}}.$$
Here the $K_i$ are nonzero constants in $k$.
If the matrix $B$ is invertible (over $\Q$),
we say that $f$ is an \textit{invertible}
Delsarte polynomial.  In this case,
we can find appropriate weights $a_0,
\ldots,a_{n+1}$ which make $f$ weighted homogeneous as follows.  Let
$q_j$ be the sum of the entries of the $j$th row of $B^{-1}$.  Then the 
equation $B B^{-1} = I_{n+2}$ means that for each $i = 0,\ldots,n+1$,
$\sum_{j = 0}^{n+1} b_{ij} q_j = 1$.  Define $d$ to be the least common
denominator of the $q_j$ and $a_j \coloneqq d q_j$.  We observe
that $d$ always divides $\det(B)$. Then $f$ is 
weighted homogeneous of degree $d$ with weights $a_0,\ldots,a_{n+1}$,
and $\gcd(a_0,\ldots,a_{n+1}) = 1$.
This collection of weights (with the gcd condition)
is uniquely determined by $f$. However, it is not 
always true that this set of weights is well-formed.

The theorem below shows that a Delsarte polynomial satisfying a 
certain gcd condition is rational.  This generalizes a result
of J. Koll\'{a}r \cite[Section 5]{Kollar} 
which proves rationality for certain hypersurfaces
defined by a loop polynomial.

\begin{theorem}
\label{thm:Delsarte}
Suppose $X \coloneqq \{f = 0\} \subset \mathbb{P}_k(a_0,\ldots,a_{n+1})$
is an 
irreducible hypersurface defined by an invertible
Delsarte polynomial in a
well-formed weighted projective space over any field $k$.  
Let $B$ be the matrix associated to $f$ and $d$ its weighted degree.
If $d = |\det(B)|$, then $X$ is rational over $k$.
\end{theorem}

\begin{proof}
Let $\mathcal{L}$ be the linear system generated by the monomials of $f$, 
where 
$$f \coloneqq \sum_{i = 0}^{n+1} K_i \prod_{j = 0}^{n+1} x_j^{b_{ij}}.$$
We claim that
this linear system induces a birational map 
$|\mathcal{L}|: \mathbb{P}(a_0,\ldots,a_{n+1}) \dashrightarrow \mathbb{P}^{n+1}$
if and only if $|\det(B)| = d$.
(Note that while $\mathcal{O}(d)$ is not necessarily a line bundle on 
$\mathbb{P}$, it is a line bundle over an open set, so a space of global
sections still induces a rational map as shown.)

To see this, view both $\mathbb{P}(a_0,\ldots,a_{n+1})$ and $\mathbb{P}^{n+1}$
as quotients of $\mathbb{A}^{n+2} \setminus \{0\}$, and let 
$M \cong \Z^{n+2}$ and
$M' \cong \Z^{n+2}$ be the lattices of the respective tori inside these
affine spaces.  The map $|\mathcal{L}|$ is induced by the map of rings
$\C[M'] \rightarrow \C[M]$ sending $y_i \mapsto \prod_{j = 0}^{n+1} 
x_j^{m_{ij}}$.  We observe that the corresponding map of lattices is 
$M' \xrightarrow{B^{\mathsf{T}}} M$.  The map on dual lattices 
is therefore $N \xrightarrow{B} N'$, and this descends to a map of quotient lattices:
$$N/(\Z \cdot (a_0,\ldots,a_{n+1})) \xrightarrow{B} N'/(\Z \cdot (1,\ldots,1)).$$
This lattice map is dual to the corresponding map of tori inside the weighted
projective spaces induced by $|\mathcal{L}|$.
It is indeed well-defined because $(a_0,\ldots,a_{n+1}) \mapsto 
(d,\ldots,d)$ and it's easy to see that the map is injective (we've used
here that $B$ is invertible). 
It is also surjective if the vectors $Be_0,\ldots,Be_{n+1},v 
\coloneqq (1,\ldots,1)$ generate $N' \cong \Z^{n+2}$, or equivalently if 
$e_0,\ldots,e_{n+1}, B^{-1} v$ generate $B^{-1} \Z^{n+2}$.  But this in turn
is the same as $B^{-1} v$ being a generator of $B^{-1} \Z^{n+2}/\Z^{n+2}$,
which is a finite abelian group of order $|\det(B)|$.  But $B^{-1} v$ is the 
vector $(a_0/d,\ldots,a_{n+1}/d)$, so its order in this group is the
least common denominator $d$.  Therefore, the map of lattices is
an isomorphism if and only if $d = |\det(B)|$.

Now suppose $d = |\det(B)|$ holds.
The transform of $X$ under $|\mathcal{L}|:\mathbb{P}(a_0,\ldots,a_{n+1}) 
\dashrightarrow \mathbb{P}^{n+1}$ is
the hyperplane $\{K_0 y_0 + \cdots + K_{n+1} y_{n+1} = 0\} \subset 
\mathbb{P}^{n+1}$.  The restriction to $X$ gives a birational map
$X \dashrightarrow \mathbb{P}^n$, completing the proof.
\end{proof}

This construction gives many new examples of rational weighted projective
hypersurfaces.  Quasismooth examples are more limited, though: a 
quasismooth Delsarte hypersurface over $\C$ has an equation which is a
sum of three types of atoms: Fermat, loop, and chain
(see \cite[Theorem 1]{KS} and \cite[Section 2.2]{ABS}).  The theorem
of Koll\'{a}r covers the case where the
defining polynomial
is a single loop, but \Cref{thm:Delsarte} gives quasismooth
examples where it is a combination of multiple loops,
for instance.  This paper focuses on finding terminal examples,
but the theorem can also be used to find new rational hypersurfaces $X$
with $K_X$ ample.

\subsection{Rational quadric bundles}
\label{subsect:quadric}

It's well-known that a smooth cubic hypersurface in 
$\mathbb{P}^{n+1}$ of even dimension
$n = 2m$ containing two disjoint $m$-planes is rational.
One way to see this is to project away from one of the planes, which gives
$X$ the birational structure of a quadric bundle over $\mathbb{P}^m$ with
a section given by the second $m$-plane.
The following
criterion can be viewed as a generalization of this fact.
Thanks to J. Koll\'{a}r for suggesting a simpler formulation of the theorem.

\begin{theorem}
\label{thm:quadric_bundle}
Let $X \coloneqq \{f = 0\} \subset \mathbb{P}_k(a_0,\ldots,a_{n+1})$ be an 
irreducible weighted projective hypersurface over a field $k$
and $1 \leq m \leq n-1$ be an integer with the following properties:
\begin{enumerate}
    \item \label{thm:quadric_bundle-1}
    $\gcd\{a_0,\ldots,a_m\} = 1$ and the set $\{a_{m+1},\ldots,a_{n+1}\}$
    is well-formed;
    \item \label{thm:quadric_bundle-2}
    Every monomial of $f$ has degree $1$ or $2$ in the variables 
    $x_0,\ldots,x_m$ and at least one has degree $1$ in these variables.
\end{enumerate}
Then $X$ is rational over $k$.
\end{theorem}

This statement also holds, suitably interpreted, in the ``degenerate" cases
$m = 0, n, n+1$.  These lead to easy rationality constructions along the lines
of \Cref{prop:lowdeg}, so we will not consider them here.

\begin{proof}[Proof of \Cref{thm:quadric_bundle}]
The idea of the proof is to look at the projection 
$X \dashrightarrow \mathbb{P}_k(a_{m+1},\ldots,a_{n+1})$.
We can understand this map better using some explicit
toric geometry. The weighted projective space
$\mathbb{P}(a_0,\ldots,a_{n+1})$ is the toric variety
corresponding to the fan which is generated by proper subsets of the 
collection of vectors $e_0,\ldots,e_{n}, v_{n+1} \in \mathbb{R}^{n+1}$ 
in the lattice $N$ generated by these vectors, where 
$$v_{n+1} \coloneqq -\frac{a_0}{a_{n+1}}e_0 - \cdots - 
\frac{a_{n}}{a_{n+1}}e_{n}.$$ 

There is a dominant toric rational map $p: \mathbb{P}(a_0,\ldots,a_{n+1}) 
\dashrightarrow \mathbb{P}(a_{m+1},\ldots,a_{n+1})$ defined 
by forgetting the first
$m+1$ coordinates.  In the toric picture, this is given by the quotient
$N \rightarrow N'$ of lattices, where 
$N' = N/(N \cap \mathrm{span}_{\R}\{e_0,\ldots,e_m\})$.
We can resolve the indeterminacy of $p$ by a single toric 
blowup $Y$ of $\mathbb{P}(a_0,\ldots,a_{n+1})$ in the stratum 
$\{x_{m+1} = \cdots = x_{n+1} = 0\}$, in particular
the one obtained by performing the barycentric subdivision of the fan in the
new ray spanned by 
$w \coloneqq -a_0 e_0 - \cdots - a_m e_m$.  Under the assumption
\eqref{thm:quadric_bundle-1}, $w$ is a primitive lattice point, and 
it's not hard to show that the quotient of lattices now gives rise to a morphism
$p: Y \rightarrow \mathbb{P}(a_{m+1},\ldots,a_{n+1})$.
\begin{commentA}

Here are the details for checking this.  The images of the rays $e_{m+1},\ldots
e_{n}, v_{n+1}$ are $\bar{e}_{m+1}, \ldots, \bar{e}_{n}, 
(-a_{m+1}\bar{e}_{m+1} - \cdots - a_{n}\bar{e}_{n+1})$ in the vector space
spanned by $e_{m+1},\ldots e_{n}, v_{n+1}$.  Notice that these are exactly
the rays used to define $\mathbb{P}(a_{m+1},\ldots,a_{n+1})$, and the lattices
match up by the assumption on well-formedness of $\{a_{m+1},\ldots,a_{n+1}\}$.

The toric rational map given by the projection of lattices is undefined only
at the strata corresponding to cones that do not map into cones of the fan of
$\mathbb{P}(a_{m+1},\ldots,a_{n+1})$.  These are the cones which contain
all of the vectors $e_{m+1},\ldots e_{n}, v_{n+1}$ among their generators.  The 
new vector $w$ is a positive linear combination of $e_{m+1},\ldots e_{n}, 
v_{n+1}$; in particular, $w = a_{m+1}e_{m+1} + \cdots + a_{n} e_{n} + 
a_{n+1} v_{n+1}$.  The new cones after the barycentric subdivision are 
precisely those that did not contain $w$ to begin with, plus for
every cone with $e_{m+1},\ldots e_{n}, v_{n+1}$ among the generators,
all the cones the
cones with one of $e_{m+1},\ldots e_{n}, v_{n+1}$ replaced by $w$.  All of
the cones of the fan defining this blowup $Y$ map into cones of the fan of
$\mathbb{P}(a_{m+1},\ldots,a_{n+1})$, proving that this blowup resolves
indeterminacies.  We note that $w$ is a primitive lattice point.  This is
because if we express $rw$ as the sum of an integer lattice point and 
$s v_{n+1}$ for $s \in \Z$, we must have $a_{n+1}|sa_i$ for $i = m+1,\ldots,
n+1$.  The assumption on well-formedness of $\{a_{m+1},\ldots,a_{n+1}\}$
gives that $s$ must be an integer.
\end{commentA}

The fibers of the toric morphism 
$p: Y \rightarrow \mathbb{P}(a_{m+1},\ldots,a_{n+1})$ can be read off from
the fan of $Y$ (for more details on how
to do this in general, see, e.g., \cite[Section 2]{HLY}).  
In our case, we only need to identify the behavior over the open torus
orbit $T \subset \mathbb{P}(a_{m+1},\ldots,a_{n+1})$. The homomorphism
of lattices associated to $p$ is surjective by 
condition \eqref{thm:quadric_bundle-1} and the collection of cones
of the fan of $Y$ contained in the preimage of $0 \in N'_{\R}$ are precisely
those generated by subsets of $e_0,\ldots,e_m, w$.  It follows that each
fiber over the open torus orbit is the toric variety corresponding to the 
fan of these cones.  Since $w = -a_0 e_0 - \cdots - a_m e_m$, this is
the weighted projective space 
$\mathbb{P}(a_0,\ldots,a_m,1)$, well-formed by \eqref{thm:quadric_bundle-1}.
In fact, $p|_{p^{-1}(T)}$ can be identified with
the second projection $\mathbb{P}(a_0,\ldots,a_m,1) 
\times T \xrightarrow{\pi_2} T$.

Now we analyze how $X$ behaves under this transformation.
If $D_i$ is the toric divisor $\{x_i = 0\}$ in 
$\mathbb{P}(a_0,\ldots,a_{n+1})$, then $D_i$ pulls back to $D_i$ in $Y$
if $i = 0,\ldots,m$ and to $D_i + a_i E$ if $i = m+1,\ldots,n+1$, where
$E$ is the exceptional divisor of the blowup.
Therefore, the total transform of $X$ in $Y$ is given by the equation
$$f(x_0,\ldots,x_m,z^{a_{m+1}}x_{m+1},\ldots,z^{a_{n+1}}x_{n+1}) = 0.$$
Here $z$ is section associated to $E$ and the left-hand side above can
be viewed as a section of the pullback of $\mathcal{O}(d)$ to $Y$.

Notice that this new equation is weighted homogeneous of degree $d$ in the 
variables $x_0,\ldots,x_m,z$ for the weighted projective space
$\mathbb{P}(a_0,\ldots,a_m,1)$, where $z$ corresponds to the weight $1$.  The 
intersection of the total transform of $X$ with each fiber over $T$ 
gives the hypersurface in $\mathbb{P}(a_0,\ldots,a_m,1)$ defined by this 
equation.
\begin{commentA}

Here are some more details for checking that this is the total transform.
A toric $\Q$-Cartier divisor $D$ corresponds to a piecewise linear function
$\varphi$ on the support of the fan,
linear on each cone. The coefficient of each boundary divisor in $D$
is determined by the value of $\varphi$ at the primitive lattice point of the
ray corresponding to $D$.  For instance, in weighted projective space
there are basic toric divisors $D_i \coloneqq \{x_i = 0\}$, 
which correspond to the 
unique piecewise function taking the value $1$ at the ray corresponding to
the divisor and $0$ on the other rays.  We can compute the pullback
of each of the divisors $D_i$ (and hence the functions $x_i$) in the blowup
$Y$ by evaluating the corresponding support functions on the new ray,
$w$.  

The ray $w$ is contained in the cone generated by 
$e_{m+1},\ldots,e_{n},v_{n+1}$, so we can see that $D_i$ pulls back to 
$D_i$ (and hence $x_i$ to $x_i$)
if $i = 0,\ldots,m$. Note that we are denoting by $D_i$ the divisor
corresponding to the same ray in $Y$.
On the other hand, $D_i$ pulls back to 
$D_i + a_i E$ and $x_i$ pulls back to $x_i z^{a_i}$ when 
$i = m+1,\ldots,n+1$, because 
$$w = a_{m+1}e_{m+1} + \cdots + a_{n+1}v_{n+1},$$
and hence
$$\varphi(w) = a_{m+1}\varphi(e_{m+1}) + \cdots + a_{n+1}\varphi(v_{n+1}).$$
Here we use $E$ to denote the exceptional toric divisor, and $z$ the 
corresponding function.
The equation for the total transform of $X$ in $Y$ is therefore
$$f(x_0,\ldots,x_m,z^{a_{m+1}}x_{m+1},\ldots,z^{a_{n+1}}x_{n+1}) = 0.$$
There will be a common factor of $z^c$ in this equation, which corresponds
to the multiple of the exceptional divisor that appears in the total 
transform.

\end{commentA}
The affine chart $\{z \neq 0\} \subset \mathbb{P}(a_0,\ldots,a_m,1)$ is
isomorphic to $\mathbb{A}^{m+1}_k$.  Restricting $p$ to this open set on each
fiber over $T$ gives
$\mathbb{A}^{m+1}_k \times T \rightarrow T$.  

The intersection $\tilde{X}$ of $X$ with $\mathbb{A}^{m+1}_k \times T$
is defined
in the fiber of $\mathbb{A}^{m+1}_k \times T \rightarrow T$
over any point
$t = [c_{m+1}: \cdots : c_{2m+1}]$ of the scheme
$T \subset \mathbb{P}
(a_{m+1},\ldots,a_{2m+1})$ (for instance the generic point) by the equation
$$\{f(x_0,\ldots,x_m,c_{m+1},\ldots,c_{2m+1}) = 0\} \subset \mathbb{A}^{m+1}_{k(t)}.$$

By condition \eqref{thm:quadric_bundle-2},
this equation is either linear or quadratic 
in the variables $x_0,\ldots,x_m$.  If it is linear, then the generic fiber
of $\tilde{X} \rightarrow T$ is isomorphic to $\mathbb{A}^m_{k(T)}$,
so $X$ is rational.  If it is quadratic, the generic fiber is an affine 
quadric over $k(T)$.  But since every monomial of $f$
has degree $1$ or $2$ in $x_0,\ldots,x_m$, the point 
$0 \in \mathbb{A}^m_{k(T)}$ is contained in this quadric, and since
at least one monomial has degree $1$ in these variables, it is smooth there.

But an affine quadric over a field $L$ containing a smooth $L$-point is 
$L$-rational, so we've shown that the
generic fiber of the affine quadric bundle $\tilde{X} \rightarrow T$ 
is rational over
$k(T)$.
The field $k(T)$ is purely transcendental over $k$, 
so it follows that $\tilde{X}$ is rational over $k$, and hence $X$ is rational
over $k$ as well.
\end{proof}

\section{Terminal examples}
\label{sect:ex}

In this section, we'll use the rationality constructions from
\Cref{sect:ratl} to prove \Cref{thm:main}.  The main sequence
of examples will come from \Cref{thm:Delsarte}.  Delsarte 
equations are of a special form, so it's not necessarily
true that a Delsarte hypersurface is very general
in its family.  Further, these examples are rarely terminal.
Nevertheless, we'll find a special choice of weights in each 
dimension $n \geq 7$ where these conditions hold.

The proposition below gives a sufficient condition for a 
hypersurface defined by a loop polynomial
(over the complex
numbers) to be terminal,
based on the Reid-Tai criterion \cite[Theorem 4.11]{Reidyoung}.
Recall that $[x]$ denotes the 
fractional part 
$x - \lfloor x \rfloor$ of $x$.

\begin{proposition}
\label{prop:reidtai}
Let $X \subset \mathbb{P}_{\mathbb{C}}(a_0,\ldots,a_{n+1})$ be a 
well-formed hypersurface defined
by a loop polynomial of type $[b_0,\ldots,b_{n+1}]$ with
$b_0,\ldots,b_{n+1} \geq 2$.  Suppose that
\begin{equation}
\label{eq:reidtaisum} \tag{II}
    \sum_{j = 2}^{n+1} [(1-b_{j-1} + b_{j-1} b_{j-2} - 
    \cdots + (-1)^{j-1} b_{j-1} b_{j-2} \cdots 
b_1)x] \geq 1 \, (\text{resp., } > 1) 
\end{equation}
for every $x \in (0,1)$.  Then $X$ is canonical (resp. terminal) in a
neighborhood of the point $[1:0:\cdots:0]$.
In particular, if this condition holds for each cyclic permutation of the
$b_j$, then $X$ is canonical (resp. terminal).
\end{proposition}

The advantage of \Cref{prop:reidtai} is that one
often needs only compute the sum of the first few terms in 
\eqref{eq:reidtaisum}.  Also, the condition does not directly
involve the weights.

\begin{proof}
A loop polynomial with exponents $b_j$ at least $2$ defines a quasismooth
hypersurface over $\C$.
Since $X$ is quasismooth at the coordinate point $[1:0:\cdots :0]$ of the 
first variable, it has a cyclic quotient singularity there. The singularity
is of type $\frac{1}{a_0}(a_2,\ldots,a_{n+1})$ \cite[Proposition 2.6]{ETW}.

We use the fact that $f$ is homogeneous of degree $d$ to give expressions
for the $a_j$, modulo $a_0$.  Since $b_0a_0 + a_1 = d$, $d \equiv a_1 
\pmod{a_0}$.  Next, $b_1 a_1 + a_2 = d$, so
$$a_2 = d - b_1 a_1 \equiv (1-b_1)a_1 \pmod{a_0}.$$
Proceeding inductively we get 
$$a_j = (1 - b_{j-1} + b_{j-1}b_{j-2} - \cdots + 
(-1)^{j-1} b_{j-1} b_{j-2} \cdots b_1)a_1\pmod{a_0}, j = 2,\ldots,n+1.$$
Note that $\gcd(a_0,a_1) = 1$, or else $a_0, a_1,$ and $d$ share a common
factor and following the loop gives that all weights share a common 
factor, a contradiction.  It follows that $a_1$ is a unit in $\mu_{a_0}$,
and multiplication by a unit does not alter the singularity type.  Hence
the singularity at the coordinate point $x_0$ is of 
type $\frac{1}{a_0}(\beta_2,\ldots,\beta_{n+1})$, where 
$$\beta_j \coloneqq 1 - b_{j-1} + b_{j-1}b_{j-2} - \cdots + 
(-1)^{j-1} b_{j-1} b_{j-2} \cdots b_1.$$
Now, the Reid-Tai criterion \cite[Theorem 4.11]{Reidyoung} states
that this singularity is canonical (resp. terminal) iff 
$$\sum_{j = 2}^{n+1} \left[\frac{i\beta_j}{a_0}\right] \geq 1  \, (\text{resp., } > 1),$$
for each $i = 1,\ldots,a_0-1$.  This sum is now the same as 
\eqref{eq:reidtaisum}, but with $x$ replaced by $i/a_0$.  Hence the 
inequality \eqref{eq:reidtaisum} for every $x \in (0,1)$
is certainly enough to imply the Reid-Tai criterion.  Finally, 
applying \cite[Proposition 2.6]{ETW} again, we note that every singularity
at a point of $X$ occurs in some stratum of the quotient singularities at
the coordinate points, so if the corresponding inequality to 
\eqref{eq:reidtaisum} holds for each coordinate point, then $X$ itself
is canonical (resp. terminal).
\end{proof}

We can now prove \Cref{thm:main}.

\begin{proof}[Proof of \Cref{thm:main}]
For each dimension $n \geq 3$, we choose $d, a_0,\ldots,a_{n+1}$ in such
a way that there is a loop polynomial
$$f \coloneqq x_0^2 x_1 + x_1^2 x_2 + \cdots + x_n^2 x_{n+1} + x_{n+1}^3 x_0,$$
which is weighted homogeneous of degree $d$.  Indeed, we can readily
compute the required weights and degree using the matrix
of the equation as in \Cref{subsect:Delsarte}.  We obtain:
\begin{align*}
    a_i & = 2^{n+1} + \sum_{j = 1}^{i} (-1)^{n+2-j} 2^{j-1},
    i = 0,\ldots,n+1, \\
    d & = 3 \cdot 2^{n+1} + (-1)^{n+2}.
\end{align*}

The first two weights differ by $1$ and $d = \det(B)$,
so the gcd condition of \Cref{thm:Delsarte} holds.

\begin{lemma}
\label{lem:monomials}
The only monomials of weighted degree $d$ in $a_0,\ldots,a_{n+1}$ are 
those in the loop polynomial $f$.
\end{lemma}

\begin{proof}
We first note that all weights but the last, $a_{n+1}$, are strictly 
between $d/4$ and $d/2$.  Indeed, the largest weight is 
$a_n < 2^{n+1} + 2^{n-1}$ so $2a_n = 2 \cdot 2^{n+1} + 2^n < d$.  
Conversely, the smallest weight besides $a_{n+1}$ is $a_{n-1}$, which
satisfies $a_{n-1} > 2^{n+1} - 2^{n-2}$, so 
$4 a_{n-1} = 4 \cdot 2^{n+1} - 2^n > d$.

This shows that any monomial of degree $d$ not involving the last variable
is the product of exactly three terms.  If $x_{k_1} x_{k_2} x_{k_3}$ has 
degree $d$ (with $0 \leq k_1 \leq k_2 \leq k_3 \leq n)$, then 
$$\sum_{j = 1}^{k_1} (-1)^{n+2-j} 2^{j-1} + 
\sum_{j = 1}^{k_2} (-1)^{n+2-j} 2^{j-1} + 
\sum_{j = 1}^{k_3} (-1)^{n+2-j} 2^{j-1} = (-1)^{n+2}.$$
We claim that $k_3 \leq k_2 + 1$.  Suppose by way of contradiction that
this is not the case.  As an alternating sum of powers of $2$, we can
see that the $k_3$ sum has absolute value greater than $2^{k_3-2}$ and 
less than $2^{k_3-1}$ (if $k_3 \geq 2)$; 
the same holds for the other sums.
Thus if $k_3 > k_2 + 1$, 
the last sum has an absolute value of more than $2^{k_3-2}$, which is 
at least twice the upper bound $2^{k_2-1}$ for the absolute value of the
other two terms.  It's therefore impossible for the sum to have absolute
value $1$ as long as $k_2 \geq 2$.  
(This leaves out the case $k_2 \leq 1$, but it's easy to check the same
is still true in that setting.)

Thus we've shown $k_2 = k_3$ or $k_2 + 1 = k_3$.  
Either way, $x_{k_1} x_{k_2} x_{k_3}$ must now be one of the loop monomials
since all the weights are distinct and the latter two terms belong to
some unique loop monomial.
Finally, we consider monomials involving the last variable $x_{n+1}$.
We note that $2a_n + a_{n+1} = d$ and $a_n$ is the largest weight, so
any monomial of degree $d$ involving $x_{n+1}$ that is not $x_n^2 x_{n+1}$
must have four terms.  We can check
$$2a_{n+1} + 2a_{n-1} \geq 
2 \cdot 2^{n+1} - 2^{n} + 2 \cdot 2^{n+1} - 2^{n-1} 
= 4 \cdot 2^{n+1} - 2^n - 2^{n-1} > d,$$
where we note that $a_{n-1}$ is the second-smallest weight.  It follows that
a monomial of degree $d$ with four terms must involve at least three copies of 
$x_{n+1}$, and hence must be the last monomial $x_{n+1}^3 x_0$.
\end{proof}

Returning to the proof of \Cref{thm:main}, \Cref{lem:monomials}
shows that every
hypersurface with degree $d$ and weights $a_0,\ldots,a_{n+1}$
has the same equation $f$, up to varying the coefficients.
The hypersurface is quasismooth if and only if all coefficients
are nonzero.
\Cref{thm:Delsarte} shows that every 
quasismooth hypersurface in this family is rational.
Every quasismooth member of the family is also a Fano variety
by the adjunction formula, since the sum of the weights is
greater than the degree.
\begin{commentA}

As an aside, the same example is quasismooth over any field, except
possibly in a few characteristics.  Indeed, the following lemma
shows that a loop polynomial of type $[b_0,\ldots,b_r]$
defines a quasismooth hypersurface over $k$ unless
$b_0 \cdots b_r = (-1)^{r+1}$ in $k$.

\begin{lemma}
\label{prop:smooth}
Let $k$ be any field.  For any positive integers $b_0,\ldots,b_r$, 
the common vanishing set of the partial derivatives of 
$x_0^{b_0}x_1 + \cdots + x_{r-1}^{b_{r-1}} x_r + x_r^{b_r} x_0$ 
in $\mathbb{A}_k^{r+1}$ is
the origin unless $b_0 \cdots b_r = (-1)^{r+1}$ in $k$.
\end{lemma}

\begin{proof}
The vanishing of the partial derivatives gives $r+1$ equations:
\begin{align*}
    & b_0 x_0^{b_0-1} x_1 + x_r^{b_r} = 0, \\
    & \hspace{3.0em} \vdots \\
    & b_{r-1} x_{r-1}^{b_{r-1}-1} x_r + x_{r-2}^{b_{r-2}} = 0, \\
    & b_r x_r^{b_r-1} x_0 + x_{r-1}^{b_{r-1}} = 0.
\end{align*}
We note that if some variable, for example $x_0$, is $0$, then the first 
equation gives that $x_r = 0$, and we can conclude by induction that
all variables vanish.  Hence any nontrivial solution to the system of 
equations above has all nonzero coordinates.  If we move the second terms 
to the right-hand side and multiply all the equations together, we obtain
$$b_0 \cdots b_r x_0^{b_0} \cdots x_r^{b_r} = (-1)^{r+1} x_0^{b_0} \cdots x_r^{b_r}.$$
Hence, if a nontrivial solution exists, we can cancel all variables 
to conclude that we must have $b_0 \cdots b_r = (-1)^{r+1}$
in $k$.
\end{proof}
    
\end{commentA}

We claim that this example
also has terminal singularities whenever 
the dimension $n \geq 7$.
This follows from \Cref{prop:reidtai} for $n \geq 8$.
Indeed, it's actually true that the sum of the first
eight terms from \eqref{eq:reidtaisum} is already greater than $1$ on the 
interval $(0,1)$
for every cyclic permutation of $[b_0,\ldots,b_{n+1}] = [2,\ldots,2,3]$,
$n \geq 8$.
Since the initial terms in the sum only depend on the first 
few of the values $b_i$, 
which are either all $2$'s or all $2$'s except a $3$ in some position, 
we only need to graph a handful of different functions to verify this.
(Or, more concretely, check a finite number of rational values in $(0,1)$
for each of these functions.)
As an example, the graph of the sum of the first five terms when 
$b_1$ through $b_5$ all equal $2$ is shown in \Cref{fig:reidtai}.
It is already greater than $1$ at every point in the open unit interval,
so we've actually shown that any coordinate point of a loop hypersurface
where the ``next five exponents" are $2$ is terminal.

\begin{figure}
    \centering
    \includegraphics[width=0.7\linewidth]{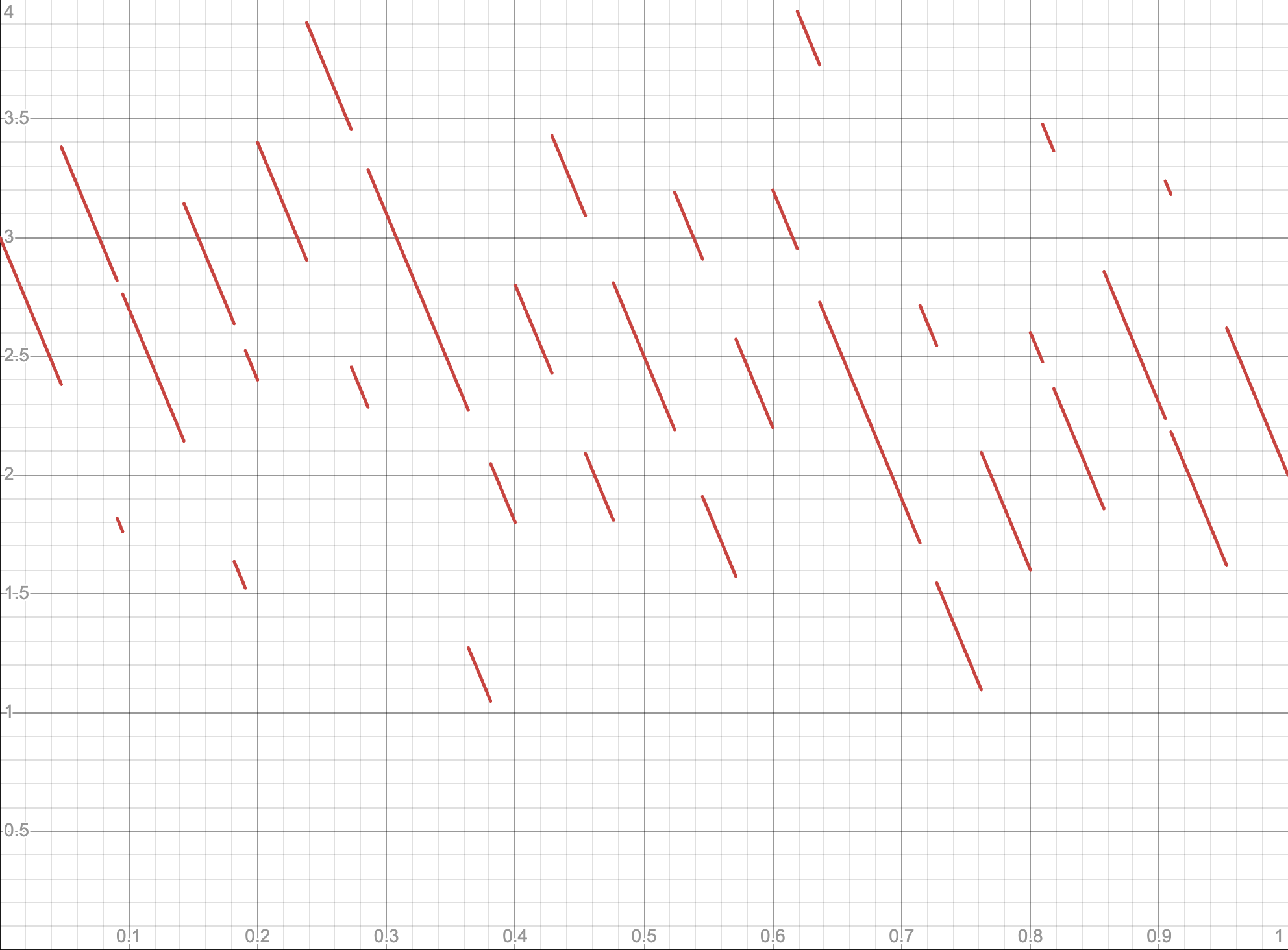}
    \caption{The sum of the first five terms in \eqref{eq:reidtaisum} for
    a sequence of exponents beginning $2,2,2,2,2,\ldots$, namely 
    the function $[-x] + [3x] + [-5x] + [11x] + [-21x]$.  This function
    is greater than $1$ on the interval $(0,1)$, 
    so the corresponding coordinate point
    is terminal.}
    \label{fig:reidtai}
\end{figure}

It turns out that the criterion \eqref{eq:reidtaisum}
fails for our example when
$n = 7$, but a computer check shows it is still terminal (using the usual
version of Reid-Tai from \Cref{thm:RT}). 
This gives an example where the criterion in \Cref{prop:reidtai} is 
not a necessary condition for $X$ to be terminal.

This completes the proof of the theorem when $n \geq 7$. 
The loop example is 
not terminal when $n = 6$, 
but the theorem is proven for $n = 6$ in \Cref{ex:non-triv_moduli}.
\end{proof}

In the loop examples in the proof above, any two quasismooth
members of the family are isomorphic to one another.  In other words,
the moduli space of quasismooth hypersurfaces is a single point.

However, we can use \Cref{thm:quadric_bundle} to give new rational examples
with non-trivial moduli answering \Cref{question:main}.  This construction
also gives the only example I could find with $n = 6$:

\begin{example}
\label{ex:non-triv_moduli}
Consider the family of hypersurfaces 
$X_{23} \subset \mathbb{P}_{\C}(9^{(2)},8^{(2)},7^{(2)},5^{(2)})$ of
dimension $6$.
In this example, the degree is more than twice the maximum
of the weights, but we claim that every quasismooth member of the family
is rational and terminal.

We first notice that every monomial of degree $23$ with the given weights
contains exactly one or two terms from the first four variables (the ones
corresponding to weights $9$ and $8$), and there exist monomials of both of
these kinds. In fact, in order for $X$ to be quasismooth, it \textit{must}
contain both monomials that are linear and 
monomials that are quadratic in the first four variables.
The general member of the family is indeed quasismooth.
Hence
the conditions of \Cref{thm:quadric_bundle} are satisfied for any
quasismooth $X$ in the family (where we take $m = 3$).
It follows that any such $X$ is rational.  A quick computation with the
Reid-Tai criterion \Cref{thm:RT} also shows that this $X$ is terminal.
This completes the proof of \Cref{thm:main} for $n = 6$.

This family 
has nontrivial moduli.  Indeed, two quasismooth hypersurfaces
in the family are isomorphic iff one is the image of the other under an
automorphism of $\mathbb{P}$ \cite[Theorem 2.1]{Esser}.
There are $26$ monomials of degree $23$,
so the projective space of weighted homogeneous polynomials
of degree $23$ is isomorphic to $\mathbb{P}^{25}$.  On the other hand,
$\mathrm{Aut}(\mathbb{P})$ has dimension $15$ (see, e.g., \cite[Lemma 1.3]{Esser}).

It's possible to find similar examples in all even dimensions $n \geq 6$.
\end{example}

\begin{remark}
\label{rem:alternative}
Another recent rationality construction due to M. Artebani and M. Chitayat 
\cite[Proposition 3.6]{Chitayat} gives an alternative approach to 
\Cref{thm:main} in the case $n \geq 7$.  Indeed, they show that 
hypersurfaces of the form $X_{ac} \subset \mathbb{P}(c^{(k)},a^{(\ell)})$,
$\gcd(a,c) = 1, k,\ell \geq 1$ are always rational, provided that the
equation $f$ involves both variables of weight $c$ and of weight $a$.
(The very general $X$ in particular has this property.)
Geometrically, their rationality proof amounts to showing that 
the projection map $X \dashrightarrow 
\mathbb{P}^{k-1} \times \mathbb{P}^{\ell-1}$ is birational.

In order to find examples failing the degree criterion 
\eqref{eq:degree_criterion}, we need $a,c > 2$, so the smallest 
possibilities are $c = 4, a = 3$.  Using the Reid-Tai criterion,
one can check that the example only has terminal singularities
when
$k > a$ and $\ell > c$.  The smallest example that
works for \Cref{thm:main} is therefore 
$X_{12} \subset \mathbb{P}(4^{(4)},3^{(5)})$, which has
dimension $n = 7$ (this is also the only example from this
construction that works
for $n = 7$). These examples have non-trivial moduli as well.
\end{remark}

\begin{commentA}
It's also possible to find rational examples
using \Cref{thm:quadric_bundle}
where the automorphism
group of the very general $X$ is trivial.

\begin{example}
\label{ex:trivauto}

Consider the following family of complex hypersurfaces of dimension $6$:
$$X_{1097} \subset \mathbb{P}^7(519,507,433,404,289,231,83,59).$$
The equation of such a hypersurface has the form
$$f \coloneqq c_0 x_0^2 x_7 + c_1 x_7^{10} x_1 + c_2 x_1^2 x_6 + c_3 x_6^8 x_2 
+ c_4 x_2^2 x_5 
+ c_5 x_5^3 x_3 + c_6 x_3^2 x_4 + c_7 x_4^2 x_0 + c_8 x_2 x_3 x_6 x_7^3 = 0.$$
The general $X$ in this family is quasismooth and the conditions of 
\Cref{thm:quadric_bundle} apply, so every quasismooth $X$ is rational.
All automorphisms of a quasismooth hypersurface in the family extend to
$\mathrm{Aut}(\mathbb{P})$ by \cite[Theorem 2.1]{Esser}.  But 
the automorphism group of this $\mathbb{P}$ is \textit{diagonal},
i.e., every element of the group is of the form $x_i \mapsto k_i x_i$ for some 
constants $k_i \in \C$.
This is because $x_i$ is the only monomial of degree $a_i$
for each $i = 0,\ldots,7$ (see also \cite[Lemma 1.3]{Esser}).
We may use the procedure of \cite[Section 7]{ETWindex} to compute the 
diagonal automorphism group of the loop hypersurface we obtain by
setting $c_8 = 0$ in the equation for $X$.  The result is that a quasismooth $X$
satisfying $c_8 = 0$ has $\mathrm{Aut}(X) \cong \mu_7$, 
where $\zeta \in \mu_7$ acts by 
$$\zeta \cdot (x_0:x_1:x_2:x_3:x_4:x_5:x_6:x_7) = (x_0:\zeta^5 x_1: \zeta^3 x_2:
\zeta^2 x_3: \zeta^4 x_4: \zeta^2 x_5: \zeta^5 x_6: \zeta x_7).$$
Indeed, we can confirm that $\zeta \cdot f = \zeta f$ under this action.
But the action on the last monomial is 
$\zeta \cdot (x_2 x_3 x_6 x_7^3) = \zeta^6 x_2 x_3 x_6 x_7^3$.  This
shows that the automorphism group of a general hypersurface in the family is
trivial.
\end{example}
\end{commentA}

We were unable to find any terminal examples answering \Cref{question:main} 
affirmatively in dimensions $4$ or $5$, using any of the constructions
mentioned above.

\end{document}